\documentclass[12pt]{amsart}
\usepackage{amscd,amssymb,verbatim}
\usepackage{amssymb,amsmath,amsthm,epsfig}

\newcommand{\N}{\mathbb{N}}

\newcommand{\R}{\mathbb{R}}

\theoremstyle{plain}
\newtheorem{Thm}{Theorem}[section]

\newtheorem{Lem}[Thm]{Lemma}
\newtheorem{Subl}[Thm]{Sublemma}

\newtheorem{Rmk}[Thm]{Remark}

\newtheorem{Def}[Thm]{Definition}
\theoremstyle{remark}

\begin{document}

\title{A note on the method of minimal vectors}
\author{George Androulakis}
\thanks{This research was partially supported by NSF}

\date{}
\maketitle

 \noindent
{\bf Abstract:} The methods of ``minimal vectors'' were introduced by
Ansari and Enflo and strengthened by Pearcy, in order to prove the existence
of hyperinvariant subspaces for certain operators on Hilbert space.
In this note we present the  method of minimal vectors for operators
on super-reflexive Banach spaces and we give a new sufficient condition
for the existence of hyperinvariant subspaces of certain operators on these
spaces..

\section{Introduction}

The {\em Invariant Subspace Problem (I.S.P.) } asks whether there exists
a separable infinite dimensional Banach space on which every operator
has a non-trivial invariant subspace. By ``operator'' we always mean
``continuous linear map'', by ``subspace'' we mean ``closed linear manifold'',
and by ``non-trivial'' we mean ``different than zero and the whole space''.
Several negative solutions to the I.S.P. are known \cite{En1} \cite{En2} \cite{R1}
\cite{R2}, \cite{R3}, \cite{R4}. It remains unknown whether the separable Hilbert
space is a positive solution to the I.S.P.. There is an extensive literature
of results towards a positive solution of the I.S.P. especially in the case
of the infinite dimensional separable complex Hilbert space $\ell_2$.
We only mention Lomonosov's result: every operator which is not a multiple of 
the identity
and  commutes with a non-zero compact operator on a complex Banach space
has a non-trivial
hyperinvariant subspace \cite{L}. For surveys on
the topic see \cite{RR} and \cite{PS}. Recently Ansari and Enflo \cite{AE}
introduced the methods of minimal vectors and gave a new proof
of the existence of non-trivial hyperinvariant subspaces
of non-zero compact operators on $\ell_2$. The method of minimal vectors which
was introduced by Enflo, was strengthened by Pearcy \cite{P} in order to
give a new proof to the following special case of Lomonosov's theorem:
every non-zero quasi-nilpotent operator on $\ell_2$ which commutes
with a non zero
compact operator has a non-trivial hyperinvariant subspace. In
this note we present the method of minimal vectors of an operator
and we carry out two
generalizations compared to the existing versions of Ansari-Enflo
and Pearcy: Firstly, the operators are  defined on a general super-reflexive
Banach space rather than the space $\ell_2$. This may be proved important
if we try to find some Banach space which is a solution to the
I.S.P. rather than examining whether $\ell_2$ is a solution to the
I.S.P.. Secondly, we introduce a property $(\star)$ that an operator may satisfy.
If an
operator $Q$ commutes with a non-zero compact operator then $Q$ satisfies
property $(\star)$. Our main result (Theorem \ref{T:main}) 
refers to operators that satisfy property $(\star)$
rather than those that commute with a non-zero compact operator.
More precisely, we prove that 
every non-zero quasi-nilpotent operator which satisfies property $(\star)$
on a super-reflexive Banach space has a non-trivial hyperinvariant subspace.
We ask whether there exist operators
which satisfy property $(\star)$ but do not have any non-zero compact operator
in their commutant. Also we ask whether every operator with no non-trivial 
invariant subspace must satisfy property $(\star)$. If the answer is positive 
then Theorem \ref{T:main} will imply that every quasi-nilpotent operator 
on a super-reflexive Banach space has a non-trivial invariant subspace. Then, 
every strictly singular operator on a super-reflexive Hereditarily
Indecomposable complex Banach space has a non-trivial invariant subspace 
(see \cite{GM}), and hence the space constructed in \cite{F} would provide 
a positive solution to the I.S.P..

We now recall some standard definitions and results that we shall use
in this paper.
A Banach space $(X, \| \cdot \|)$ is called {\em strictly convex\/} if for 
every $x, y \in X$ with $\| x \| = \| y \| = \| (x + y ) / 2 \| =1$ we have that 
$x=y$. 
A Banach space $(X, \|\cdot\|)$ is called {\em uniformly
  convex\/} if for every $\varepsilon>0$
there exists a $\delta=\delta(\varepsilon)>0$ such that
for $x,y\in X$ with $\| x \| = \|y\|=1$ and $\left\|\frac{x+y}2\right\| > 1- \delta$
we have that $\|x-y\| < \varepsilon$. The function $\delta(\varepsilon)$
is called the modulus of uniform convexity of $X$.
The norm of $X$ is {\em G\^ateaux differentiable} if
for every $x\in X\backslash\{0\}$ and for every $y\in X$ the limit
\begin{equation}\label{eq1}
\lim_{t\to 0} \frac{\|x+ty\| - \|x\|}t
\end{equation}
exists. The Banach space $X$ is called
{\em smooth\/} if for every $x\in X\backslash\{0\}$
there exists a unique $f \in X^*$ such that $f(x) = \|x\|^2
= \| f \|^2$. We denote the functional $f$ by $(x)^*$.
It can be proved that the norm of $X$ is  G\^ateaux
differentiable if and only if $X$ is smooth, in which case  the limit
in (\ref{eq1}) is equal to $\text{Re}\,(x)^*(y)$.
The norm of $X$ is called {\em Fr\'{e}chet differentiable\/} if the
limit in (\ref{eq1}) exists uniformly for all $y \in X$ with
$\| y \| =1$. The norm of $X$ is called {\em uniformly smooth\/}  if the 
limit in (\ref{eq1}) exists uniformly for all $x,y \in X$ with 
$\| x \| = \| y \| =1$. 
A Banach space $X$ is called {\em super-reflexive\/} if every infinite
dimensional space $Y$ which is finitely represented in $X$ must
be reflexive. It is proved in \cite{En0} (see also \cite{Pi}) 
that every super-reflexive Banach space $X$ can be equivalently
renormed to be uniformly convex. It follows from
a renorming technique of Asplund \cite{A}  that a Banach space is
super-reflexive if and only if it can be equivalently renormed to be
uniformly convex or uniformly smooth or both.

\section{Minimal vectors and invariant subspaces}\label{sec2}

We start by introducing some notations and terminology. If $X$ is a
Banach space, $x \in X$ and $\varepsilon >0$ we denote by
$\text{S}(x, \varepsilon)$ (respectively $\text{Ba}(x, \varepsilon)$) the {\em sphere\/}
(respectively the {\em closed ball\/}) of $X$ with center $x$ and radius
$\varepsilon$, namely the set $\{ y \in X: \| x-y \| = \varepsilon \}$
(respectively the set $\{ y \in X: \| x - y \| \le \varepsilon \}$).

\begin{Def} Let $X$ be a Banach space and $Q$ be an operator on $X$.
  We say that an operator $Q$ satisfies property $(\star)$ if for every
  $\varepsilon \in (0,1)$ there exists $x_0 \in \text{S}(0,1)$ such
  that for every weakly convergent sequence
$(x_n) \subset \text{S}(x_0, \varepsilon)$ there
  exists a subsequence $(x_{n_k})_k$ of $(x_n)$ and a sequence $(K_k)
  \subset \{ Q \}'$ such that
\begin{itemize}
\item[(a)] $\| K_k \| \le 1$ and $K_k(x_0) \ge \frac{1+ \varepsilon}2$
  for all $k \in \N$.
\item[(b)] $( K_k(x_{n_k}) )_k$ converges in norm.
\end{itemize}
\end{Def}

The purpose of (a) is to ensure that the limit of (b) is non-zero. Indeed,
notice that if $Q, \varepsilon, x_0, (K_k)_k, (x_{n_k})_k$ are as in
the previous definition then 
$\| K_kx_0 -K_kx_{n_k} \| \le \| K_k \| \| x_0-x_{n_k} \| \le \varepsilon$,
thus
$\| K_k x_{n_k} \| \ge \| K_k x_0 \| - \| K_kx_0 -K_kx_{n_k}\|
\ge (1+ \varepsilon)/2 - \varepsilon = (1- \varepsilon) /2 >0$.

Also notice that for every operator $Q$ if there exists a non-zero
compact operator which commutes with $Q$ then $Q$ satisfies property
$(\star)$.

Our main result is the following:

\begin{Thm} \label{T:main}
  Let $X$ be a super-reflexive Banach space and $Q$ be a non-zero
  quasi-nilpotent operator on $X$ which satisfies property $(\star)$. Then $Q$ has
  a non-trivial hyperinvariant subspace.
\end{Thm}

For the proof of this result we use the method of minimal
  vectors of an operator.  Let $(X, \|\cdot\|)$ be a reflexive
Banach space, $Q$ be an operator on $X$ with dense range, $x_0 \in X$
with $\| x_0 \|=1$ and $\varepsilon \in (0,1)$.  We define a
sequence of {\em minimal vectors of $Q$ with
respect to $x_0$ and $\varepsilon$} to be a sequence $(y_n)_{n\in\mathbb N}$
as follows: For every $n \in \N$
the set $Q^{-n}\text{Ba}(x_0, \varepsilon)$ is non-empty (since $Q^n$
has a dense range), closed and convex. Thus there exists $y_n \in
Q^{-n}\text{Ba}(x_0, \varepsilon)$ such that
\begin{equation} \label{eq2}
\| y_n \| = \inf \{ \| y \| :
y \in Q^{-n}\text{Ba}(x_0, \varepsilon) \} .
\end{equation}
Indeed, if $(y_{n,m})_m$ is a sequence in $Q^{-n}\text{Ba}(x_0,
\varepsilon)$ with
\begin{equation} \label{eq3}
\| y_{n,m} \| \searrow \inf \{ \| y \| :
y \in Q^{-n}\text{Ba}(x_0, \varepsilon) \} ,
\end{equation}
then $(y_{n,m})_m$ is a subset of $Q^{-n}\text{Ba}(x_0,
\varepsilon)\cap \text{Ba}(0, \| y_{n,1}\|)$ which is weakly compact
(since it is a closed, convex and bounded subset of a reflexive
space).  Thus by passing to a subsequence and relabeling we can assume
that $(y_{n,m})_m$ converges weakly to some vector $y_n$. Since the
norm is weakly lower semicontinuous, (\ref{eq3}) implies (\ref{eq2}).

In order to prove Theorem \ref{T:main} we need the following three
results whose proofs are postponed.
For the first result, notice that if $X$ is a reflexive Banach space,
$Q$ is an operator on $X$ with dense range, $x_0 \in X$ with $\| x_0 \|=1$,
$\varepsilon \in (0,1)$ and $(y_n)$ is a sequence of minimal vectors of
$Q$ with respect to $x_0$ and $\varepsilon$, then
the sequence $((Q^ny_n-x_0)^*)_n$ is bounded, (namely,
$\|(Q^ny_n-x_0)^*\| = \|Q^ny_n-x_0\| = \varepsilon)$
by the minimality of $\|y_n\|$), thus it has
weak$^*$ limit points. We want to know that 0 is not a weak$^*$ limit
point of the sequence $((Q^ny_n-x_0)^*)_n$. The next result
yields that this is true provided that the choice of $\varepsilon$ is appropriate.

\begin{Lem} \label{L:unifconvex}
  Let $(X, \|\cdot\|)$ be a smooth and uniformly convex Banach
  space and $Q$ be an operator on $X$ with dense range.  
Then there exists $\varepsilon\in\left[\frac12, 1\right)$ such that 
the following is satisfied: if $x_0 \in X$
  with $\| x_0 \|=1$, $(y_n)_n$ is a sequence of minimal
  vectors of $Q$ with respect to $x_0$ and $\varepsilon$, and $f$ is a
  weak$^*$ limit point of $(( Q^ny_n-x_0)^*)_n$, then $f\ne 0$.
\end{Lem}

\begin{Lem}\label{L:quasinilpotent}
  Let $X$ be a reflexive Banach space, $Q$ be quasi-nilpotent operator
  on $X$ with dense range, $x_0 \in X$ with $\| x_0 \|=1$, $\varepsilon
  >0$, and $(y_n)_{n\in \mathbb N}$ be a sequence of minimal vectors
  of $Q$ with respect to $x_0$ and $\varepsilon$.  Then there exists
  an increasing sequence $(n_k)_k$ of $\N$ such
  that
\begin{equation}\label{eq4}
\lim_k \frac{\|y_{n_k-1}\|}{\|y_{n_k}\|} = 0.
\end{equation}
\end{Lem}

For the next Lemma, if $X$ is a Banach space and $f\in X^*$ then
$\ker(f)$ denotes the {\em kernel\/} of $f$.

\begin{Lem}\label{L:perpendicular}
  Let $X$ be a reflexive smooth Banach space, $Q$ be an operator on
  $X$ with dense range, $x_0 \in X$ with $\| x_0 \|=1$, $\varepsilon
  \in (0,1)$, and $(y_n)_{n\in\mathbb N}$ be a sequence of minimal
  vectors of $Q$ with respect to $x_0$ and $\varepsilon$.  Then for
  $n\in\mathbb N$,
\begin{equation}\label{eq5}
\ker( ( y_n )^*) \subseteq \ker((Q^n)^* (Q^ny_n-x_0 )^*).
\end{equation}
\end{Lem}

Now we are ready for the

\begin{proof}[Proof of Theorem \ref{T:main}]
  Since $X$ is super-reflexive we can assume by our discussion in the
  previous section, that $(X, \| \cdot \|)$ is smooth and locally
  uniformly convex.  Without loss of generality we assume that $Q$ has
  a dense range and it is 1-1 (because the range and the kernel of $Q$
  are hyperinvariant subspaces of $Q$).  By Lemma~\ref{L:unifconvex}
there exists $\varepsilon \in \left[\frac{1}{2},1\right)$ such that the conclusion of the
lemma is satisfied. For that $\varepsilon$, since $Q$ satisfies property $(\star)$, 
let $x_0\in X$, $\|x_0\| = 1$
  such that the statement of the definition of property $(\star)$ is valid
for the operator $Q$. Let $(y_n)_n$ be a sequence of
  minimal vectors of $Q$ with respect to $x_0$ and
  $\varepsilon$. By Lemma~\ref{L:quasinilpotent} let $(n_k)_k$ be an increasing
  subsequence of $\N$ such that (\ref{eq4}) is valid. Since $X$ is reflexive,
by considering a further subsequence of $(n_k)$ and relabeling we can
assume that $(Q^{n_k -1}y_{n_k-1})_k$ converges weakly. By the property $(\star)$
of $Q$, there exists a subsequence of $(n_k)$ (which, by relabeling, is
still called $(n_k)$) and a sequence $(K_k)_k \subset \{ Q \}'$ such that
$(K_k Q^{n_k-1}y_{n_k -1})_k$ converges in norm to some vector
$w \in X$.  By our discussion following the definition of property $(\star)$
we have that $w \not = 0$. Since $Q$ is 1-1 we have that $Qw \not = 0$.
We claim that $Y:= \{ Q \}' (Qw)$ is a non-trivial hyperinvariant subspace
for $Q$. We only need to show that $Y \not =X$. For that reason we let
$f$ to be a weak$^*$ limit point of $((Q^{n_k}y_{n_k} -x_0)^*)_k$,
which is non-zero by Lemma~\ref{L:unifconvex}, and we will show
that $Y \subset \ker (f)$. We need to show that if $T \in \{ Q \}'$ then
$f(TQw)=0$. Let $T \in \{ Q \}'$ and $k \in \N$.
Since $\ker((y_{n_k})^*)$ is a 1-codimensional subspace of $X$ and
  $y_{n_k} \notin \ker((y_{n_k})^*)$ (notice that
$(y_{n_k})^* (y_{n_k}) = \|y_{n_k}\|^2\ne 0$), we have that
$X = \text{span}\{y_{n_k}\} \oplus \ker((y_{n_k})^*)$, thus there
  exists a scalar $a_k$ and $r_k\in \ker((y_{n_k})^*)$ such that
\begin{equation}\label{eq6}
TK_k (y_{n_k-1}) = a_k y_{n_k} + r_k.
\end{equation}
We claim that $a_k \to 0$. Indeed,
\begin{align*}
  |a_k| \|y_{n_k}\|^2 &= | (y_{n_k})^* (a_ky_{n_k} + r_k)|\\
  &= |(y_{n_k})^* T K_k (y_{n_k-1})| \quad \text{(by (\ref{eq6}))}\\
  &\le \| (y_{n_k})^*\| \| T \| \| K_k \| \|y_{n_k-1}\|\\
  &\le \|y_{n_k}\| \|T \| \|y_{n_k-1}\|,
\end{align*}
thus  $|a_k| \le \| T \| \|y_{n_k-1}\|/\|y_{n_k}\|
\mathop{\longrightarrow} \limits_{k\to\infty} 0$.

First apply $Q^{n_k}$ and then $( Q^{n_k}y_{n_k} - x_0)^*$,
on (\ref{eq6}), to obtain
\begin{gather}
Q^{n_k} T K_k y_{n_k-1} =a_k Q^{n_k}y_{n_k} + Q^{n_k} r_k \quad \text{and}\nonumber\\
\label{eq7}
(Q^{n_k}y_{n_k}-x_0)^* Q^{n_k} T K_k y_{n_k-1} =
a_k (Q^{n_k}y_{n_k}-x_0)^* Q^{n_k}y_{n_k} +
(Q^{n_k}y_{n_k}-x_0)^* Q^{n_k} r_k .
\end{gather}
Since $r_k \in\ker((y_{n_k})^*)$, by Lemma~\ref{L:perpendicular}
we have that $(Q^{n_k}y_{n_k} - x_0)^* Q^{n_k}(r_k) = 0$.
Furthermore, since $a_k \to 0$ and
$\|Q^{n_k}y_{n_k} - x_0\| = \varepsilon$, we have that the right hand side of
(\ref{eq7}) tends
to zero.  Thus by taking limits and noticing that $T,K_k \in \{ Q\}'$
for all $k$, (\ref{eq7}) becomes
\begin{equation}\label{eq8}
\lim_k (Q^{n_k}y_{n_k} - x_0)^* TQ K_k Q^{n_k-1} y_{n_k-1} = 0.
\end{equation}
Since $(K_kQ^{n_k-1}y_{n_k-1})_k$ converges in norm to $w$ and
$f$ is a weak$^*$ limit point of $((Q^{n_k}y_{n_k} -x_0)^*)_k$,
(\ref{eq8}) yields that $f(TQw)=0$ which finishes the proof.
\end{proof}

We now turn our attention to the proof of Lemma~\ref{L:unifconvex}.
Before presenting the proof of Lemma~\ref{L:unifconvex} we need
the following two results.

\begin{Subl}\label{slem8}
  Let $(X, \|\cdot\|)$ be a smooth and strictly convex Banach space,
  $x_0\in X$, $\|x_0\| = 1$, $0<\varepsilon<1$ and $w\in
  \text{S}(x_0,\varepsilon)$. The following conditions are equivalent:
\begin{itemize}
\item[(a)] $\|x_0-\lambda w\| > \varepsilon$ for all
 $\lambda\in [0,1)$.
\item[(b)] $\text{Re}\,\frac{(x_0-w)^*}{\|x_0-w\|^2} (x_0) \ge 1$.
\end{itemize}
\end{Subl}

\begin{proof}
  (a) $\Rightarrow$ (b):\ Since $\|x_0-\lambda w\| > \varepsilon =
  \|x_0-1\cdot w\|$, for all $\lambda\in [0,1)$, we have that the
  derivative of the function $f(\lambda) = \|x_0- \lambda w\|$ at 1 is
  non-positive. Set $g(\mu) = f(1-\mu)$. Then $g'(0) = -f'(1)$, thus
  $g'(0)\ge 0$. Note that $g(\mu) = \|x_0-(1-\mu)w\| = \|x_0-w + \mu
  w\|$. Thus
\begin{align*}
  g'(0) &= \text{Re}\, \frac{(x_0-w)^*}{\| (x_0-w)^*\|} (w) = \text{Re}\,\frac{(x_0-w)^*}{\|x_0-w\|}(w)\\
  &= \text{Re}\, \frac{(x_0-w)^*}{\|x_0-w\|} (w-x_0+x_0) = -\frac{\|x_0-w\|^2}{\|x_0-w\|} + \text{Re}\, \frac{(x_0-w)^*}{\|x_0-w\|}(x_0)\\
  &= -\|x_0-w\| + \text{Re}\, \frac{(x_0-w)^*}{\|x_0-w\|} (x_0).
\end{align*}
Thus $g'(0) \ge 0$ if and only if
\[
\text{Re}\, \frac{(x_0-w)^*}{\|x_0-w\|^2} (x_0) \ge 1.
\]
(b) $\Rightarrow$ (a):\ By the proof of
(a)~$\Rightarrow$~(b), notice that if (b) is valid then $f'(1) \le 0$.
Notice also that $f(1) = \varepsilon$ since $w\in \text{S}(x_0,\varepsilon)$.
Thus if (b) is valid then $f(\lambda)\ge \varepsilon$ for
$\lambda\in [\lambda_0,1)$ for some $\lambda_0<1$.  On the other hand,
there are at most two numbers $\lambda_1\le \lambda_2$ such that
$f(\lambda_1) = f(\lambda_2) = \varepsilon$ and
$f(\lambda) < \varepsilon$ for $\lambda$ between
$\lambda_1,\lambda_2$ (since $(X, \|\cdot\|$) is strictly convex).
Since $f(\lambda) \ge \varepsilon$ for $\lambda\in [\lambda_0,1)$, it
is impossible to have $\lambda_1 < \lambda_2 = 1$.  Hence
$f(\lambda) \ne \varepsilon$ for $\lambda<1$. Since
$f(0) = \|x_0\| = 1 > \varepsilon$, we have that
$f(\lambda) > \varepsilon$ for $\lambda\in [0,1)$.
\end{proof}

\begin{Lem}\label{lem7}
  If $(X, \|\cdot\|)$ is a smooth and uniformly convex Banach space.
  Then for every $\eta>0$ there
  exists $\varepsilon\in \left[\frac12,1 \right)$ such that for all
$x_0 \in X$ with $\| x_0 \| =1$ and for all  
$w\in \text{S}(x_0,\varepsilon)$ satisfying $\|x_0-\lambda w\|
  >\varepsilon$ for all $\lambda\in [0,1)$, we have that
$\|w\| \le  \eta$.
\end{Lem}

\begin{proof}
  We start with the

\noindent
{\bf Claim:}\ Let $(X, \|\cdot\|)$ is a uniformly
convex Banach space, and $\eta'>0$.
Let $\delta( \cdot )$ denote the modulus of uniform convexity of $X$. 
Then for every $x_0 \in X$ with $\| x_0 \|=1$,
\begin{equation}\label{eq9}
\sup\{\|x\|\colon \ x\in \text{S}(x_0,1),  \ \text{Re}\,(x_0-x)^* (x_0)
> 1-2\delta(\eta' ) \} < \eta'.
\end{equation}
Indeed, if $x \in \text{S}(x_0,1)$  with
$\| x \| \ge \eta'$, then
we have that $\|x_0-(x_0-x)\| \ge \eta', \|x_0\| = 1, \|x_0-x\| = 1$,
hence
\begin{align*}
1- \delta(\eta') & \ge
  \left\|\frac{x_0+(x_0-x)}2\right\|
  \ge \left|(x_0 - x)^* \left(\frac{x_0+(x_0-x)}2\right)\right|\\
  &\ge \text{Re}\,(x_0-x)^* \left(\frac{x_0+(x_0-x)}2\right)\\
  &= \frac{\text{Re}\,(x_0-x)^*(x_0) + 1}2.
\end{align*}
Thus $\text{Re}\,(x_0-x)^*(x_0) \le 1-2 \delta(\eta')$ which
finishes the proof of the Claim.

Let $X$ and $\eta$ as in the statement of Lemma~\ref{lem7}.
Let $\varepsilon$ satisfying
\begin{equation}\label{eq10}
 \frac{1}{2} \le \varepsilon, \quad
0< \frac1\varepsilon - 1 \le 2\delta \left(\frac{\eta}{2}\right) \quad \text{and}
 \quad 1-\varepsilon \le \eta/2.
\end{equation}
Let $w\in \text{S}(x_0,\varepsilon)$ satisfying
$\|x_0-\lambda w\| > \varepsilon$ for all $\lambda\in [0,1)$.
By Sublemma~\ref{slem8} we
have that
\begin{equation}\label{eq11}
\text{Re}\,\frac{(x_0-w)^*}{\|x_0-w\|^2} (x_0) \ge 1.
\end{equation}
Let $x\in \text{S}(x_0,1)$ with
\begin{equation}\label{eq12}
x_0-w = \varepsilon(x_0-x).
\end{equation}
Then
\begin{align*}
 \left|\text{Re}\,\frac{(x_0-w) ^*}{\|x_0-w\|^2} (x_0) -
   \text{Re}\,(x_0-x)^* (x_0)\right|
 &\le \left|\frac{(x_0-w)^*}{\|x_0-w\|^2} (x_0) - (x_0-x)^* (x_0)\right|\\
  &\le \|x_0\| \left\|\frac{(x_0-w)^*}{\|x_0-w\|^2} - (x_0-x)^*\right\|\\
 &= \left\|\frac{\varepsilon(x_0-x)^*}{\|\varepsilon(x_0-x)\|^2} -
   (x_0-x)^*\right\| \quad \text{(by (\ref{eq12}))}\\
 &\le \left| \frac{1}{\varepsilon} -1 \right| \le
 2 \delta \left(\frac{\eta}{2} \right)  \quad \text{(by (\ref{eq10}))}.
\end{align*}
Therefore by (\ref{eq11}) we have that
$\text{Re}\,(x_0-x)^* (x_0)\ge 1-2 \delta \left( \frac{\eta}{2} \right)$. 
By (\ref{eq9})
we have that $\|x\| < \eta/2$. By (\ref{eq12}) we have that
$w =(1-\varepsilon) x_0 + \varepsilon x$ and thus by the triangle
inequality, $\| x \| < \eta/2$, and (\ref{eq10}) we obtain
$\|w\| \le (1-\varepsilon) + \varepsilon\|x\| < (1-\varepsilon) + \|x\| < \eta$.
\end{proof}

 \begin{proof}[Proof of Lemma \ref{L:unifconvex}]
  Let $X$ be a smooth and 
  uniformly convex Banach space and $Q$ be an
  operator on $X$ with dense range. For $\eta = \frac13$ we choose
  $\varepsilon\in\left[\frac12, 1\right)$ to satisfy the statement of
  Lemma~\ref{lem7}.

  Let $(y_n)_{n\in\mathbb N}$ be a sequence of minimal vectors of
  $Q$.  For $n\in\mathbb N$, by the minimality of $\| y_n \|$ we have that 
 $Q^ny_n \in \text{S}(x_0,\varepsilon)$
  and $\|x_0-\lambda Q^ny_n\|>\varepsilon$ for all
   $\lambda\in[0,1)$. Thus by Lemma~\ref{lem7} we obtain
   $\|Q^ny_n\|\le \frac13$.
  Let $f$ be any weak$^*$ limit point of the sequence
  $((Q^ny_n-x_0)^*)_{n\in \mathbb N}$. Since $\|x_0\|=1$,
  $\|Q^ny_n\|\le \frac13$ and $\|Q^ny_n-x_0\| = \varepsilon \ge 1/2$ we
  have that $f\ne 0$. Indeed for $n\in \mathbb N$,
\begin{align*}
  (Q^ny_n-x_0)^*(-x_0) &= (Q^ny_n-x_0)^* (Q^ny_n-x_0 - Q^ny_n)\\
  &= ( Q^ny_n - x_0)^* (Q^ny_n-x_0) - ( Q^ny_n-x_0)^* (Q^ny_n)\\
  &= \|Q^ny_n-x_0\|^2 - (Q^ny_n-x_0)^*(Q^ny_n)\\
  &\ge \varepsilon^2 - |(Q^ny_n-x_0)^* (Q^ny_n)|\\
  &\ge \varepsilon^2 - \|(Q^ny_n-x_0)^*\| \|Q^ny_n\|\\
  &\ge \varepsilon^2 - \varepsilon\cdot\frac13\\
  &\ge \frac14 - \frac16 = \frac1{12} > 0 \quad \left(\text{since }
    \frac12 \le \varepsilon\right).
\end{align*}
Since $f$ is a weak$^*$ limit point of $((Q^ny_n-x_0)^*_n$, we have
that $f(-x_0) \ge \frac1{12}$, thus $f\ne 0$.
\end{proof}

\begin{proof}[Proof of Lemma \ref{L:quasinilpotent}]
  If the statement is not true, then there exists a positive number
  $\delta$ such that
\[
\frac{\|y_{n-1}\|}{\|y_n\|} \ge \delta \text{ for all } n \in \N .
\]
Thus for every integer $n\in \mathbb N$, we have
\begin{equation}\label{eq13}
\|y_1\| \ge \delta \|y_2\| \ge \delta^2\|y_3\|\ge \cdots\ge \delta^n \|y_{n+1}\|.
\end{equation}
We have that $\|Qy_1 - x_0\|\le \varepsilon$. Also we have that
$\|Q(Q^ny_{n+1}) - x_0\| = \|Q^{n+1} y_{n+1} - x_0\|\le \varepsilon$.
By the minimality of $\|y_1\|$ we have:
\begin{equation}\label{eq14}
\|y_1\| \le \|Q^n y_{n+1}\| \le \|Q^n\| \|y_{n+1}\|.
\end{equation}
By combining (\ref{eq13}) and (\ref{eq14}) we have $\delta \le
\|Q^n\|^{1/n}$ which is a contradiction since $Q$ is quasi-nilpotent.
\end{proof}

In order to prove Lemma \ref{L:perpendicular} we need the following

\begin{Rmk} \label{R:exercise}
  Let $X$ be a Banach space, $f\in X^*\backslash \{ 0 \}$ and $g \in X^*$
such that
\begin{equation}\label{eq15}
\text{for all } x\in X, \text{ if \rm Re}(f(x)) < 0 \text{ then \rm Re}(g(x))\ge 0.
\end{equation}
Then $\ker(f) \subseteq \ker(g)$ and moreover there exists a non-positive 
real number  $a$ such that $g=af$.
\end{Rmk}

\begin{proof}
Let $x_0 \in X$ with $f(x_0)=-1$. Thus $\text{Re}\,g(x_0) \ge 0$. 
We first claim that $\ker(f) \subseteq \ker(g)$. Indeed, otherwise there exists 
$x \in \ker(f) \backslash \ker(g)$. Without loss of generality assume that
$\text{Re}\,g(x) < -2 \text{Re}\,g(x_0)$. Let $x'= x_0 + x$. Then 
$f(x')= f(x_0) + f(x)= f(x_0)$, hence $\text{Re}\,f(x') =-1 <0$. We also have
$\text{Re}\, g(x')= \text{Re}\, g(x_0)+ \text{Re}\, g(x) < - \text{Re}\, g(x_0) \le 0$ 
which is a contradiction, proving that $\ker(f) \subseteq \ker(g)$.

Since both $\ker(f), \ker(g)$ are at most $1$-codimensional subspaces 
of $X$ we have that there exists a scalar $a$ such that $g=a f$. 
Since $f(x_0)=-1$ and $\text{RE}\, g(x_0) \ge 0$ we have that
$\text{Re}\, a \le 0$. If $a \not \in \R$ then let $a= a_1 + i a_2$ with 
$a_1, a_2 \in \R$, $a_1 \le 0$ and $a_2 \not = 0$. Let $x_1 \in X$
with $f(x_1)= -1 + ia_2^{-1}(1-a_1)$. Then $\text{Re}\, f(x_1)= -1 <0$
and $\text{Re}\, g(x_1)= \text{Re}\,((a_1 +i a_2)(-1 + i a_2^{-1}(1-a_1)))=-1<0$
which is a contradiction, proving that $a \in \R$. 
\end{proof}

Now we are ready for the

\begin{proof}[Proof of Lemma \ref{L:perpendicular}]
For a fixed $n \in \N$ we prove that the assumption of Remark~\ref{R:exercise}
is satisfied for $f=(y_n)^*$ and $g=(Q^n)^* ((Q^ny_n -x_0)^*)$. Let
$x \in X$ with $\text{Re}\,(y_n)^*(x) < 0$. We claim that 
$\text{Re}\, (Q^n)^*(Q^ny_n-x_0)^* (x) \ge 0$ i.e. 
$\text{Re}\, (Q^ny_y-x_0)^*(Q^nx) \ge 0$. 
 Indeed, otherwise, since
$\text{Re}\,(Q^ny_n-x_0)^* (Q^nx)$ is
the derivative of the function
\[
t \mapsto \|Q^ny_n-x_0 + tQ^nx\|
\]
at 0, we obtain that this function is decreasing for $t$ in a
neighborhood of 0. Thus for small $t>0$ we have
\[
\varepsilon = \|Q^ny_n-x_0\|\ge \|Q^ny_n-x_0 + tQ^nx\|
\]
i.e.
\[
\|Q^n(y_n+tx) -x_0\| \le \varepsilon.
\]
We have by the
minimality of $\|y_n\|$ that
\[
\|y_n\| \le \|y_n + tx\| \quad \text{for small}\quad t>0.
\]
Thus the derivative of the function
\[
t\mapsto \|y_n + tx\|
\]
must be non-negative at 0, i.e.\ $\text{Re}\,(y_n)^* (x) \ge 0$ which
is a contradiction.
\end{proof}

 \vspace{.2in} \scriptsize{
\noindent
Department of Mathematics, University of South Carolina, Columbia, SC
29208.  giorgis@math.sc.edu

\end{document}